RESEARCH(Research Manuscript)

# Improved Multi-Dimensional Bee Colony Algorithm for Airport Freight Station Scheduling

Haiquan Wang[1,*], Menghao Su[2], Ran Zhao[1], Xiaobin Xu[3], Hans-Dietrich Haasis[4], Jianhua Wei[2], Shengjun Wen[1], Yan Wang[2], Ping Liu[1] and Hongjun Li[1]

## Abstract

Due to the rapid increase of air cargo and postal transport volume, an efficient automated multi-dimensional warehouse with elevating transfer vehicles (ETVs) should be established and an effective scheduling strategy should be designed for improving the cargo handling efficiency. In this paper, artificial bee colony algorithm, which possesses strong global optimization ability and fewer parameters, is firstly introduced to simultaneously optimize the route of ETV and the assignment of entrances and exits. Moreover, for further improve the optimization performance of ABC, novel full-dimensional search strategy with parallelization, and random multi-dimensional search strategy are incorporated in the framework of ABC to improve the diversity of the population and the convergence speed respectively. Our proposed algorithms are evaluated on several benchmark functions, and then applied to solve the combinatorial optimization problem with multitask, multiple entrances and exits in air cargo terminal. The simulations show that the proposed algorithms can achieve much more desired performance than the traditional artificial bee colony algorithm at balancing the exploitation and exploration abilities.

## Keywords

Scheduling; Artificial bee colony algorithm; Air cargo terminal; Parallelization; Random multi-dimensional search

# 1. Introduction

Nowadays, air transport activities around the world develop rapidly, it is estimated that half the value of all international trade will be moved by air within the next decade. To improve the cargo processing speed and avoid congestion problem, the management of air cargo terminal which is the major gateway for cargo services become crucial, and an effective and efficient scheduling strategy considering the route of ETVs and the assignment of entrances and exits should be designed to minimize the time cost for handling all inbound and outbound cargoes [1-3].

For the scheduling problem in air cargo terminal with multiple entrances and exits which is a NP-hard problem, corresponding research have received considerable attention across the past few years. Guo [4] and Qiu [5] studied the inbound and outbound cargoes scheduling problem with single ETV and solved it with genetic algorithm and particle swarm optimization (PSO) algorithm respectively. PSO was also applied to assign two ETVs to different cargo areas with improved shared fitness strategy in [6]. Lei [7] discussed the task scheduling problem with ETVs considering collision avoidance and introduced expert system to improve the operation efficiency.

The works mentioned above focus on deciding the cargoes' transportation sequence with considering picking sequence and ETV routing, but very few research discussed the problem of assignment of entrances and exits if there are several gates in the airport freight station. In our research, a swarm intelligent algorithm named artificial bee colony algorithm (ABC), which possesses strong global optimization ability and fewer parameters [8-11] is firstly proposed and adopted to solve the problem of scheduling the entrances and exits with several outbound and inbound tasks and the action of ETV simultaneously.

*Corresponding Author: Haiquan Wang (wanghq@zut.edu.cn)
[1] Zhongyuan Petersburg Aviation College, Zhongyuan University of Technology, 41 Zhongyuan Road, Zhengzhou 450007, China
[2] Faculty of Electrical and Engineering, Zhongyuan University of Technology, 41 Zhongyuan Road, Zhengzhou, 450007, China
[3] School of Automation, Hangzhou Dianzi University, 115 Wenyi Road, Xihu District, Hangzhou, China
[4] Business Studies and Economics, University of Bremen, Zinckestraße 20, 28865 Lilienthal, Germany

Actually, ABC has been well adapted for various complex optimization problems, such as hotel recommendation [12], civil engineering design [13], aerospace industry [14], software testing [15], logistics warehouse management [16], manufacturing production [17], communication problem [18] and so on. However, it often suffers from the problem of premature convergence because of its single-dimensional random search strategy in the employed bee phase and the onlooker bee phase. To accelerate the convergence speed meanwhile obtaining high accuracy solutions, other metaheuristic algorithms were introduced and combined with traditional ABC to improve its performance. Ustun & Toktas [19] replaced the onlooker bee operator with the mutation and crossover phases of differential evolution algorithm to increase the accuracy and speed up the convergence. Aiming at improving the convergence accuracy of ABC algorithm, combined with the learning characteristics of Q-learning algorithm, the update dimensions in each iteration of ABC algorithm could be dynamically adjusted in [20]. Xu et al. [21] introduced differential evolution strategy in employed bee phase to accelerate its convergence and adopted the global best position to guide the following updating processes in onlooker bee phase which could enhance the local search ability. The firework explosion search mechanism was introduced to explore the potential food sources of ABC in [22]. A modified search operator was employed in [23] to exploit useful information of current best solution in onlooker phase, and the ability of exploitation could be improved. Obviously, most of them are based on the introduced search strategies or the combination with other algorithms, and there are few systematic analyses and improvements from the perspective of operation mechanism. Therefore, for balancing the abilities of exploration and exploitation, after modeling the actions of ETV in cargo terminal with multiple entrances and exits, improved ABC algorithms with paralleled full dimensional search strategy and random multi-dimensional search are proposed in this paper.

The rest of paper is organized as follows: Section 2 introduces the established scheduling model of air cargo terminal. Then two improved strategies based on ABC algorithm are proposed and verified in Section 3. In Section 4, the improved algorithms are applied to solve the ETV scheduling problem considering multiple entrances and exits. Finally, the above work is summarized.

# 2. The scheduling model of the airport freight station

The airport freight station consists of container storage area, bulk cargo storage area, and unhandled cargo area. As the core of the whole station, the container storage area is used for handling the containerized cargoes, which are unloaded from aircraft in the airside or inbounded from the bulk cargo storage area in the landside. It is a three-dimensional warehouse with 16 entrances and exits and two rows of shelves, each row has eight layers and 60 columns, the total slots are $60 \times 8 \times 2 = 960$.

ETV is employed for handling cargoes among different entrances and exits. Fig. 1 shows the operation procedure of ETV if there are three inbound and outbound tasks needed to be delivered. Firstly, ETV receives the assignment and moves to the entrance $R_1$, picks up the good and places it at $I_1$. And then ETV determines the next outbound task and loads the good at $O_1$, sends it to the scheduled exit $C_1$. After finishing the current outbound task, ETV picks up the task $O_2$ and carries it to the exit $C_2$.

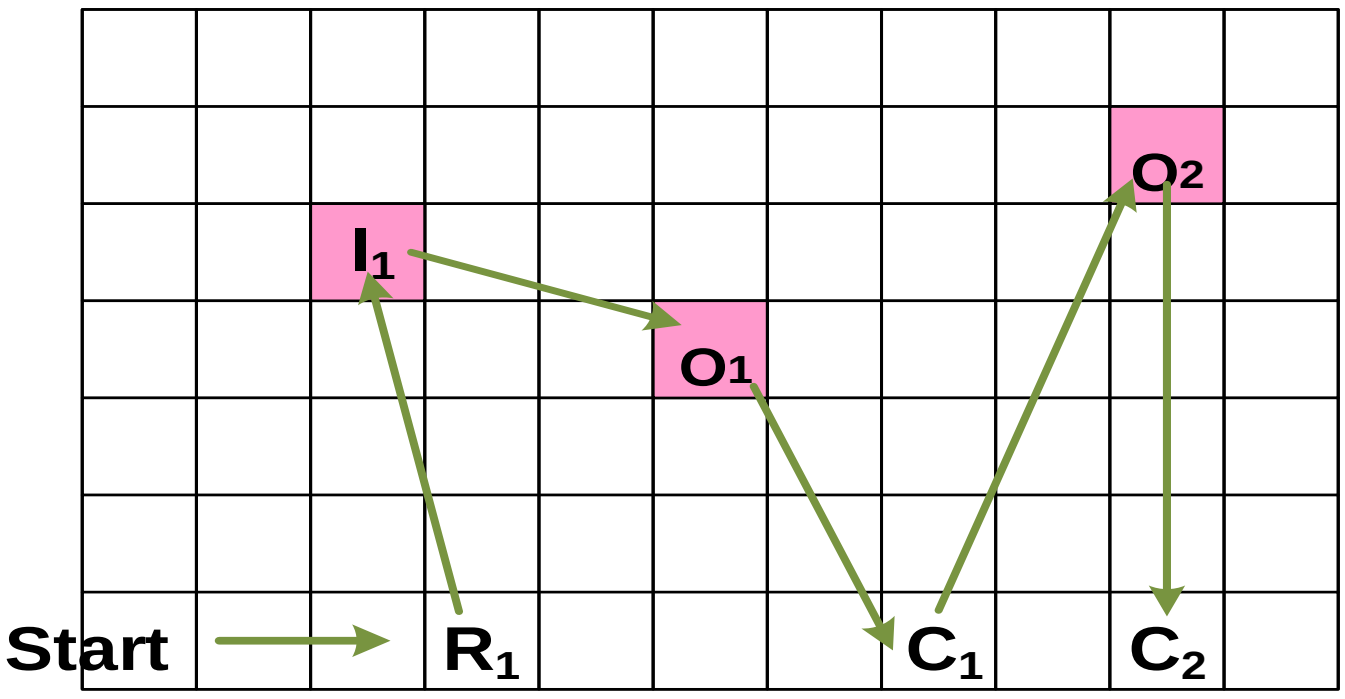

**Fig. 1.** The operation procedure of ETV.

During the pickup and delivery operation, ETV always goes through three phases, they are acceleration, constant speed, and deceleration in horizontal and vertical directions respectively. Thus the time needed to finish a whole pickup and delivery process is determined by the maximum one between the horizontal travelling time $T_x$ and the vertical lifting time $T_y$, they are defined as follows

$$\begin{cases} T_x = 2 * \sqrt{\frac{w*u}{a_x}}, \; w * u \le D_x; \\ T_x = T_{xi} + \frac{w*u-D_x}{V_x}, \; w * u > D_x. \end{cases} \quad (1)$$

$$\begin{cases} T_y = 2 * \sqrt{\frac{h*e}{a_y}}, h * e \le D_y; \\ T_y = T_{yj} + \frac{h*e-D_y}{V_y}, h * e > D_y. \end{cases} \quad (2)$$

Where $e$ and $u$ are the differences of layer and column respectively between the starting and ending points, $h$ and $w$ are the height and width of storage location, $a_x$ and $a_y$, $V_x$ and $V_y$ are the accelerations and speeds in horizontal and vertical directions. $T_{xi}$ and $T_{yj}$, $D_x$ and $D_y$ are corresponding time costs and the travelled distances when ETV accelerates to maximum speed and immediately decelerates to 0 in different directions, the values of them could be determine with the following formulas.

$$T_{xi} = 2 * V_x/a_x \;, \; T_{yj} = 2 * V_y/a_y \quad (3)$$

$$D_x = \frac{1}{4} * \; a_x * \; {T_{xi}}^2, D_y = \frac{1}{4} * \; a_y * \; {T_{yj}}^2 \quad (4)$$

With the motion analysis of ETV, the time costs from any position to another could be calculated and part of them are listed in Table 1 (Only the values corresponding to first five layers and six columns are listed). The row or column of the matrix represents the number of differences of rows or columns between the starting and ending positions.

**Table 1.** The matrix of time costs.

| | 1 | 2 | 3 | 4 | 5 | 6 |
|---|---|---|---|---|---|---|
| 1 | 0 | 5.47 | 7.74 | 9.62 | 11.50 | 13.37 |
| 2 | 11.62 | 11.62 | 11.62 | 11.62 | 11.62 | 13.37 |
| 3 | 22.87 | 22.87 | 22.87 | 22.87 | 22.87 | 22.87 |
| 4 | 34.12 | 34.12 | 34.12 | 34.12 | 34.12 | 34.12 |
| 5 | 45.37 | 45.37 | 45.37 | 45.37 | 45.37 | 45.37 |

Thus, the time cost $J_i$ to execute the $i^{th}$ task, which is the total time cost for picking up, releasing as well as transporting all assigned cargoes, is expressed as,

$$J_i = J_{i_0} + J_{i_1} + 2\zeta \quad (5)$$

Here $J_{i_0}$ is the time needed to arrive at cargo's position from its current position, $J_{i_1}$ is the time to move to the destination position after ETV get the cargo, ξ is the time for ETV to load or unload the cargo. $J_{i_0}$ and $J_{i_1}$ can be obtained from the matrix of time cost of ETV in Table 1.

Based on the above analysis, if there are several inbound and outbound tasks in the cargo terminal, the key of improving the efficiency of cargo transportation is to schedule the sequence of inbound/outbound tasks as well as the action of ETV, and reduce the total time cost as defined below.

$$\sum_{i=1}^{n} J_i = \sum_{i=1}^{n} (J_{i_0} + J_{i_1}) + 2n\zeta \quad (6)$$

Obviously, it is a typical optimization problem where the total time cost in (6) is the fitness function, and an effective optimization algorithm should be introduced to solve it.

# 3. Methodologies

## 3.1 ABC algorithm

As a kind of swarm intelligence algorithm, ABC algorithm simulates the actions of natural bees where possible solution is represented by food source, the quality of associated solution is equal to the nectar amount of food source.

The algorithm begins with a randomly distributed initial population generation and evaluation [8]. Then repeated search cycles execute to update the population of solutions. During the cycles, the employed bee probabilistically produces a neighbor food source $x'_{ik}$ around current optimal solution $x_{ik}$ with (7), and the fitness values $fitness_m$ of generated solutions will be evaluated to determine if the new one is better than the original solution.

$$x'_{ik} = x_{ik} + rand(-1,1) * \left(x_{ik} - x_{jk}\right) \quad (7)$$

Where $i$ and $j$ represent the number of the solutions, $i, j \in \{1,2...,N\}$, $i \neq j$. $k$ is the dimension of solution which is selected randomly, $k \in \{1,2..., D\}$.

The $m^{th}$ onlooker bee randomly chooses to exploit or not around corresponding employed bee with the probability $P_m$ defined as follow

$$P_m = \frac{fitness_m}{\sum_{m=1}^{SN} fitness_m} \quad (8)$$

If current solution to be exploited cannot improve for several iterations, it will be abandoned, and a randomly produced scout bee will replace it.

## 3.2 The improved ABC algorithm

Obviously, in ABC algorithm, a random single dimension search is executed which means only one dimension is randomly selected and updated according to (7) in employed bee phase and onlooker bee phase. In this case, the updated dimension may be different in each iteration and the optimal dimension obtained in the previous iterations is likely to be omitted in the following iterations. Thus, the search toward the possible optimal solution is unable to be continued, the optimization accuracy and the convergence speed will be affected. To improve the diversity of the population and the probability of obtaining the optimal solution, meanwhile improving the convergence speed, several improvements are introduced in this paper.

### 3.2.1 Paralleled full-dimensional ABC algorithm

A full-dimensional artificial bee colony (fdABC) algorithm is introduced in this paper, different from traditional ABC, it will traverse all dimensions of the solution and select the optimal dimension for further exploration, therefore the search could be extended and the possibility of obtaining optimal solution will be improved. The main steps of fdABC algorithm are shown as follows.

Step 1: Initialization: Initialize the parameters of ABC, such as the food source, the population and the maximum number of iterations.

Step 2: Employed bee phase: Execute the neighborhood search as (7) where *k* varies from zero to the number of dimensions. If the fitness value of the updated solution is better than the previous one, the solution with the updated dimension will be retained and further search will be done based on the obtained one until all dimensions are searched.

Step 3: Onlooker bee phase: Onlooker bee selects food source with roulette strategy and execute full-dimensional search around the selected solution as step 2.

Step 4: Scout bee phase: If the iteration reaches the threshold of limit and no better solution is found, the scout bees can update new solutions randomly.

Step 5: The global optimal solution obtained so far will be recorded and the algorithm jumps to step 2 for further exploration until reaching the maximum iteration number.

Because more dimensions need to be explored during the updating processes, fdABC could improve the probability of obtaining the optimal solution and avoid premature convergence. Conversely, as it should travel more dimensions, the time cost will increase inevitably.

For further improving its efficiency, a master-slave parallel mode [24] is applied to the most time-consuming stages such as the initial fitness function calculation and the updating procedure in employed bee phase, in that mode the population is divided into different parts and then the repeated calculations are executed in different cores of multicore processor. Therefore, the corresponding tasks could be finished in parallel meanwhile the structure of fdABC algorithm will not be affected. The necessary steps of parallel full-dimensional ABC algorithm (PfdABC) are arranged as follows:

Step 1: Initialization.

Step 2: The initial population is divided into six parts, and the fitness values of each solution are calculated in different cores of the CPU.

Step 3: Updating in employed bee phase. All the employed bees are equally distributed into the cores of CPU and full-dimensional searches are performed for each group in different threads.

The update strategies of the onlooker bee phase and scout bee phase are the same as traditional ABC algorithm.

### 3.2.2 Random multi-dimensional artificial bee colony algorithm

fdABC algorithm travels all the dimensions for each updating process, which could improve the optimized accuracy but affect its efficiency, another algorithm named random multi-dimensional artificial bee colony (RmdABC) algorithm is proposed to balance the abilities of exploration and exploitation. The key improvement of the strategy is to randomly select several different dimensions from {1,2,…,D} for one solution, and execute the updated process with (7) in the employed bee and onlooker bee phases. The number of updated cycles for discovering a new food source depends on the amounts selected from the set of {1,2,…,D}, it could be one to D. Obviously, RmdABC randomly traverses any several dimensions of the solution in each iteration. On the one hand, fewer dimensions are updated compared with fdABC, and its time complexity could be greatly improved. On the other hand, it covers more dimensions for each solution compared with ABC, which means it possesses higher possibility to obtain the optimal solution. The pseudo-code of RmdABC is shown as follow:

**Table 2.** The pseudo-code of RmdABC.

| RmdABC |
|---|
| 01: //Initialization, |

```
set the maximum number of iterations, the swarm size N, the number of dimension D
    for i = 1 to FoodNumber
      flag = 0;
      Random_D = rand_i(D);      //Generate a random sequence
        //Random multi-dimensional greedy search strategy
        for j = 1 to Random_D
            produce the candidate solution with (7) and evaluate its fitness value;
            if fitness (Sol_i) < fitness (Food_i) then Food_i = Sol_i and flag = 1;
        end for
        if flag = 1 then trial = 0; else trial +1;
   end for
    According to (8), calculate the probability prob_i to determine if onlooker bee chooses to exploit
    or not around the ith employed bee
   for i = 1 to FoodNumber
      flag = 0;
      if rand < prob_i
          produce the candidate solution with (7) and evaluate its fitness value;
        if fitness (Sol_i) < fitness (Food_i) then Food_i = Sol_i and flag = 1;
        if flag = 1 then trial = 0; else trial +1;
     end if
   end for
   //Scout bee phase
   if trial > Limit
     trial = 0;
     Randomly generate a solution;
   end if
end for
```

As stated above, more dimensions will be explored in fdABC compared with ABC, but more time consuming is needed. With the master-slave parallel strategy, calculations could be separated into different cores of CPU, the time cost for optimization could be reduced; through random dimensional selection, less dimensions will be explored which means the efficiency will be improved. Therefore, the proposed two strategies could effectively balance the abilities of exploration and exploitation.

## 3.3 Performance verification of the proposed algorithms

To verify the performance of the improved ABC algorithms, benchmark functions as shown in Table 3 are introduced, and ABC, fdABC, PfdABC as well as RmdABC are applied to solve them. The parameters of algorithms are set as follows: The swarm size is set to be 200, the maximum number of local searches is 100, and the maximum number of iterations is set as being 1000. The algorithms are executed ten times, and the simulation results, including average runtime, shortest runtime, average optimal results, best optimization solutions and their variances are shown in Tables 4-6 corresponding to different dimensions.

**Table 3.** Benchmark functions

| The name of function | Formulae | Search range | Optimal value |
| --- | --- | --- | --- |
| Bent Cigar Function | $f_1(x) = x_1^2 + 10^6 \sum_{i=2}^{D} x_i^2$ | $[-100,100]^D$ | 0 |
| Sum of Different Power Function | $f_2(x) = \sum_{i=1}^{D} \lvert x_i \rvert^{i+1}$ | $[-100,100]^D$ | 0 |
| Rosenbrock's Function | $f_3(x) = \sum_{i=1}^{D-1} (100({x_i}^2 - x_{i+1})^2 + (x_i - 1)^2)$ | $[-100,100]^D$ | 0 |
| Rastrigin's Function | $f_4(x) = \sum_{i=1}^{n} [x_i^2 - 10\cos(2\pi x_i) + 10]$ | $[-500,500]^D$ | 0 |
| Step's Function | $f_5(x) = \sum_{i=1}^{n} (\lvert x_i + 0.5 \rvert)^2$ | $[-100,100]^D$ | 0 |

**Table 4.** The performance of PfdABC, RfdABC, fdABC and ABC with 60 dimensions.

| The name of function | Algorithm | Average Runtime/s | Average Optimal Result | Best Optimal Result | Shortest Runtime/s | Var of Optimal Result |
|---|---|---|---|---|---|---|
| Bent Cigar Function | ABC | **111.849** | 2.320e+5 | 1.709e+5 | **111.4955** | 2.039e+9 |
| | fdABC | 253.316 | **2.929e-252** | **1.416e-252** | 251.577 | **0** |
| | RmdABC | 188.172 | 1.654e-2 | 1.49e-4 | 179.199 | 2.28e-4 |
| | PfdABC | 256.408 | 2.995e-16 | 1.028e-16 | 179.199 | 1.055e-32 |
| Sum of Different Power Function n | ABC | **121.258** | 3.024e+41 | 1.593e+36 | **120.91** | 4.433e+83 |
| | fdABC | 676.282 | **8.680e-255** | **7.047e-258** | 672.29 | **0** |
| | RmdABC | 401.083 | 2.091e-42 | 6.660e-46 | 395.373 | 1.791e-83 |
| | PfdABC | 490.647 | 4.164e-89 | 1.902e-90 | 488.905 | 4.053e-177 |
| Rosenbrock's Function | ABC | **116.400** | 5490.448 | 4884.468 | **115.536** | 2.083e+5 |
| | fdABC | 284.881 | **2.33e-3** | **1.27e-4** | 279.176 | **2.021e-06** |
| | RmdABC | 212.479 | 1.1359 | 0.323 | 208.025 | 1.013 |
| | PfdABC | 288.203 | 0.140 | 0.0072 | 286.3592 | 0.0245 |
| Rastrigin's Function | ABC | **121.356** | 198.642 | 179.002 | **121.051** | 163.148 |
| | fdABC | 250.787 | **0** | **0** | 245.425 | **0** |
| | RmdABC | 190.324 | 4.145e-06 | 2.271e-07 | 189.681 | 2.759e-11 |
| | PfdABC | 279.681 | 2.956e-13 | 2.274e-13 | 272.530 | 3.447e-27 |
| Step's Function | ABC | **120.309** | 0.232 | 0.177 | **119.949** | 0.843e-4 |
| | fdABC | 224.032 | **0** | **0** | 218.470 | **0** |
| | RmdABC | 178.819 | 1.262e-08 | 4.285e-10 | 178.200 | 1.121e-16 |
| | PfdABC | 266.942 | 3.298e-22 | 1.464e-22 | 255.218 | 1.416e-44 |

**Table 5.** The performance of PfdABC, RfdABC, fdABC and ABC with 80 dimensions.

| The name of function | Algorithm | Average Runtime/s | Average Optimal Result | Best Optimal Result | Shortest Runtime/s | Var Optimal Result |
|---|---|---|---|---|---|---|
| Bent Cigar Function | ABC | **186.356** | 7.443e+7 | 6.612e+7 | **185.785** | 8.560e+13 |
| | fdABC | 349.896 | **9.500e-25** | **1.073e-251** | 341.500 | **0** |
| | RmdABC | 254.843 | 1.830 | 0.187 | 242.967 | 8.657926 |
| | PfdABC | 336.140 | 3.739e-11 | 2.729e-11 | 330.117 | 6.778e-23 |
| Sum of Different Power Function | ABC | **161.169** | 7.666e+82 | 9.065e+73 | **157.306** | 5.875e+166 |
| | fdABC | 1092.133 | **5.377e-250** | **9.849e-253** | 1083.415 | **0** |
| | RmdABC | 641.927 | 1.236e-18 | 4.733e-22 | 629.479 | 1.193259e-35 |
| | PfdABC | 767.770 | 9.062e-72 | 2.578e-72 | 761.159 | 6.124e-143 |
| Rosenbrock's Function | ABC | **153.868** | 65508.192 | 53885.588 | **152.849** | 7.237e+7 |
| | fdABC | 408.311 | **0.00517** | **0.000291** | 393.921 | **3.638e-05** |
| | RmdABC | 297.845 | 2.626 | 1.1049 | 287.474 | 1.758 |
| | PfdABC | 388.867 | 0.264 | 0.041 | 380.135 | 0.064 |
| Rastrigin's Function | ABC | **155.747** | 759.215 | 687.486 | **153.785** | 2578.573 |
| | fdABC | 351.982 | **0** | **0** | 346.306 | **0** |
| | RmdABC | 275.060 | 0.0184 | 0.0056 | 270.149 | 9.127e-05 |
| | PfdABC | 358.622 | 4.206e-13 | 2.274e-13 | 353.856 | 5.888e-27 |
| Step's Function | ABC | **151.843** | 8.142 | 6.349 | **151.444** | 0.965 |
| | fdABC | 300.048 | **0** | **0** | 297.340 | **0** |
| | RmdABC | 537.590 | 6.049e-07 | 2.1807e-08 | 373.611 | 3.606e-13 |
| | PfdABC | 333.573 | 3.940e-17 | 2.531e-17 | 323.896 | 9.741e-35 |

**Table 6.** The performance of PfdABC, RfdABC, fdABC and ABC with 100 dimensions.

| The name of function | Algorithm | Average Runtime/s | Average Optimal Result | Best Optimal Result | Shortest Runtime/s | Var Optimal Result |
|---|---|---|---|---|---|---|
| 1.Bent Cigar Function | ABC | **223.196** | 3.787e+8 | 3.416e+8 | **222.774** | 1.099e+15 |
| | fdABC | 425.692 | **2.044e-251** | **1.091e-251** | 420.020 | **0** |
| | RmdABC | 303.243 | 12.037 | 0.405 | 300.169 | 199.478 |
| | PfdABC | 399.895 | 4.392e-08 | 3.463e-08 | 387.853 | 6.669e-17 |
| 2.Sum of Different Power Function | ABC | **199.808** | 2.279e+120 | 7.215e+111 | **197.343** | 2.222e+241 |
| | fdABC | 1664.371 | **1.671e-245±0** | **6.863e-248** | 1619.398 | **0** |
| | RmdABC | 920.823 | 3.470e-06 | 4.151e-09 | 911.149 | 1.069e-10 |
| | PfdABC | 1120.629 | 8.680e-57 | 3.9115e-58 | 1115.067 | 5.118e-113 |
| 3.Rosenbrock's Function | ABC | **189.997** | 5.919e+3 | 4.346e+5 | **186.093** | 9.399e+9 |
| | fdABC | 538.976 | **0.00884** | **0.001** | 504.637 | **1.04e-4** |
| | RmdABC | 351.792 | 6.042937 | 1.430 | 337.071 | 18.386 |
| | PfdABC | 456.292 | 1.021 | 0.201 | 449.671 | 1.568 |
| 4.Rastrigin's Function | ABC | **193.990** | 152.642 | 193.287 | **193.2873** | 65.193 |
| | fdABC | 458.470157 | **0** | **0** | 452.714 | **0** |
| | RmdABC | 574.398 | 7.97e-4 | 8.599e-05 | 573.802 | 7.859e-07 |
| | PfdABC | 454.618 | 1.136e-10 | 8.356e-11 | 450.936 | 2.190e-22 |
| 5. Step's Function | ABC | **189.754** | 78.395 | 63.6545 | **187.726** | 97.986519 |
| | fdABC | 637.867 | **0** | **0** | 622.008 | **0** |
| | RmdABC | 289.625 | 1.379e-07 | 1.838e-08 | 287.601 | 1.87107e-14 |
| | PfdABC | 398.953 | 5.632e-14 | 2.467e-14 | 391.252 | 3.979e-28 |

Table 4, Table 5 and Table 6 illustrate the optimization results with four different algorithms in different dimensions, where the best results among the five indices are highlighted in bold font. As can be seen from the tables, all algorithms could solve the nonlinear problems within limited time and the corresponding optimal variances on different dimensions show the stability of ABC could be greatly improved with the other three strategies. Furthermore, an increase in problem dimension did not require exponential increment in population size or evaluation number, it can be stated that the proposed algorithms are not very sensitive to increments in problem dimensions and has a good scalability.

As seen from the tables, the average runtime and the shortest runtime increase from 60 dimensions to 100 dimensions because more dimensions need to be updated for all algorithms. Especially for fdABC, the execution time for solving Step's function increases by 203% with respect to growing dimensions.

For the three indices including average optimal result, best optimal result and the variant of result, the values corresponding to fdABC is smallest for all functions with different dimensions, but its time complexity are the highest for most functions, thus the ability of exploration could be proved. For the Step's and Rastrigin's Functions in lower dimensions, PfdABC needs longer runtime than fdABC because the communications between different cores of CPU takes much more complex than calculations in lower dimensions.

Another improved algorithm RmdABC gets smaller fitness values and variant values but longer time compared with the ones of ABC. On the other hand, it improves the time complexity of fdABC at the expense of a worse solution because it updates less dimensions.

Fig. 2 graphically shows the curves of fitness values for Expanded Schaffer's Function with different dimensions. Obviously, fdABC, PfdABC and RmdABC algorithms have faster convergence speed and better solutions compared with ABC.

Based on the results above, RmdABC, PfdABC, fdABC strategies improves the performance and the stability of ABC effectively. fdABC possesses the best optimization result, but it spends too much time on optimizing. The optimization solutions corresponding to PfdABC and RmdABC are not as good as fdABC, but they can reduce the optimization time greatly. The time cost of RmdABC is shorter than fdABC and PfdABC because RmdABC randomly selects dimensions to be updated in each iteration. Obviously, the proposed PfdABC and RmdABC algorithms could balance the abilities of exploration and exploitation compared with ABC and fdABC, and the corresponding algorithms could be selected based on different requirements with considering the time cost and accuracy.

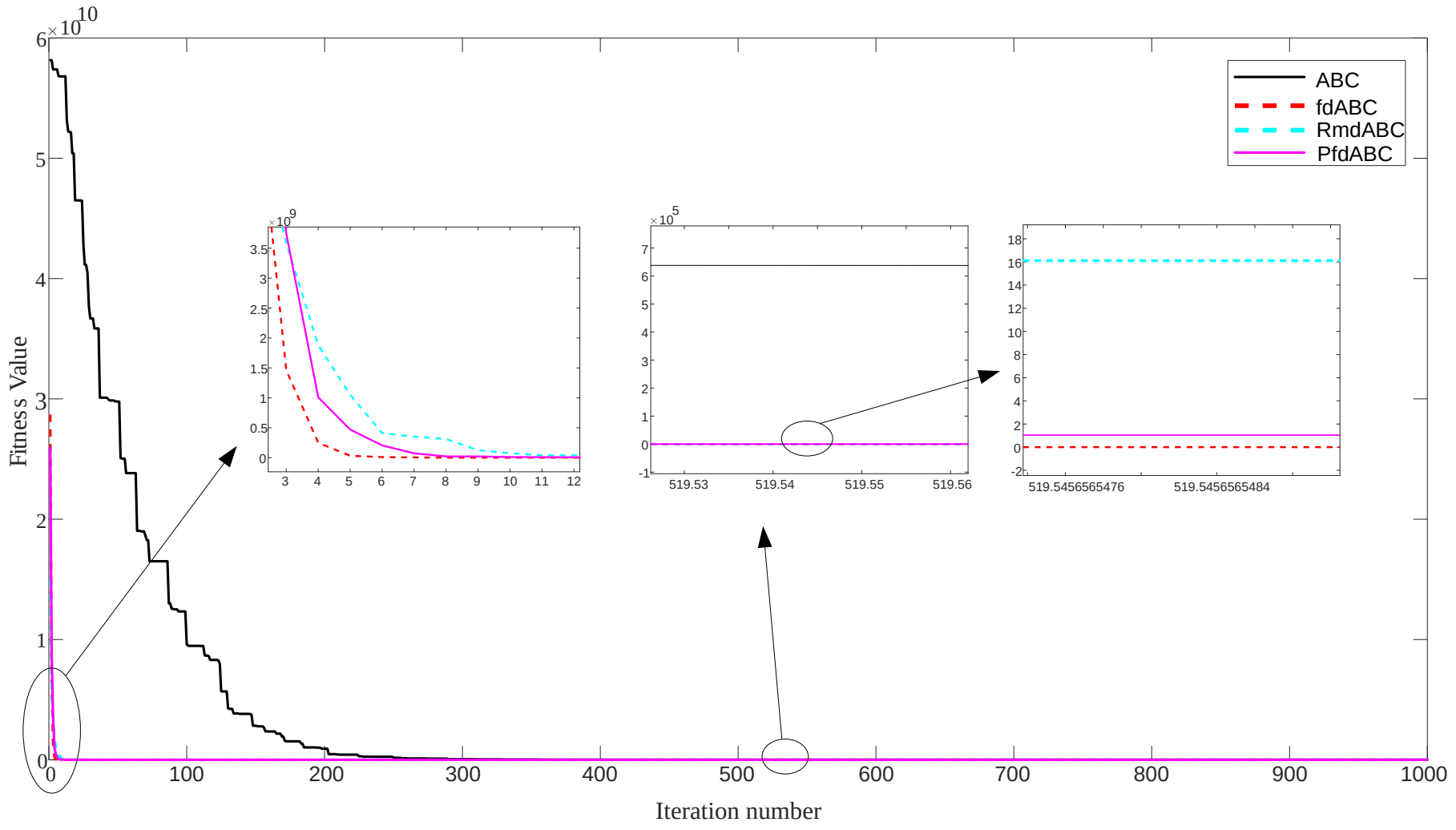

(60)

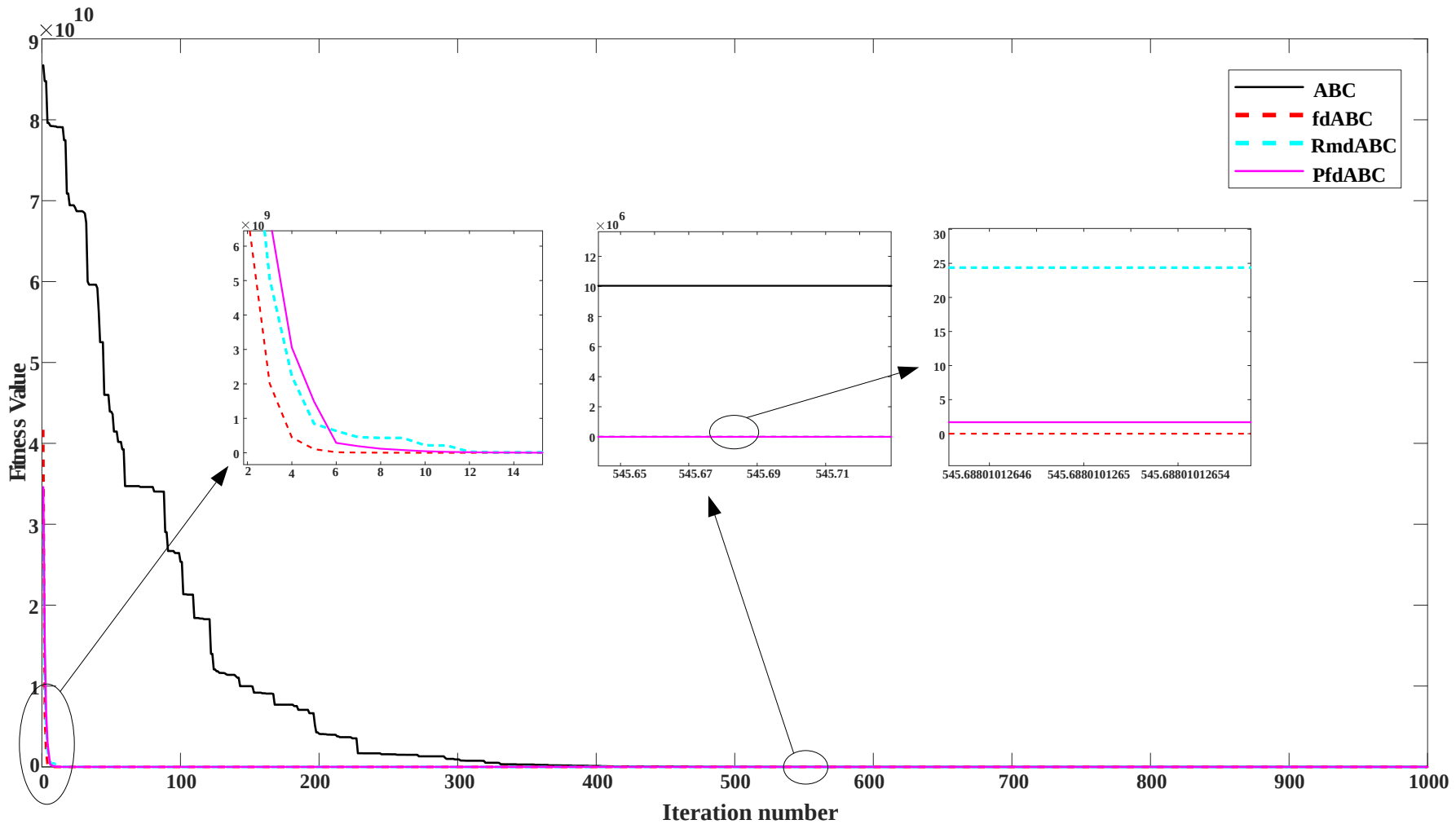

(80)

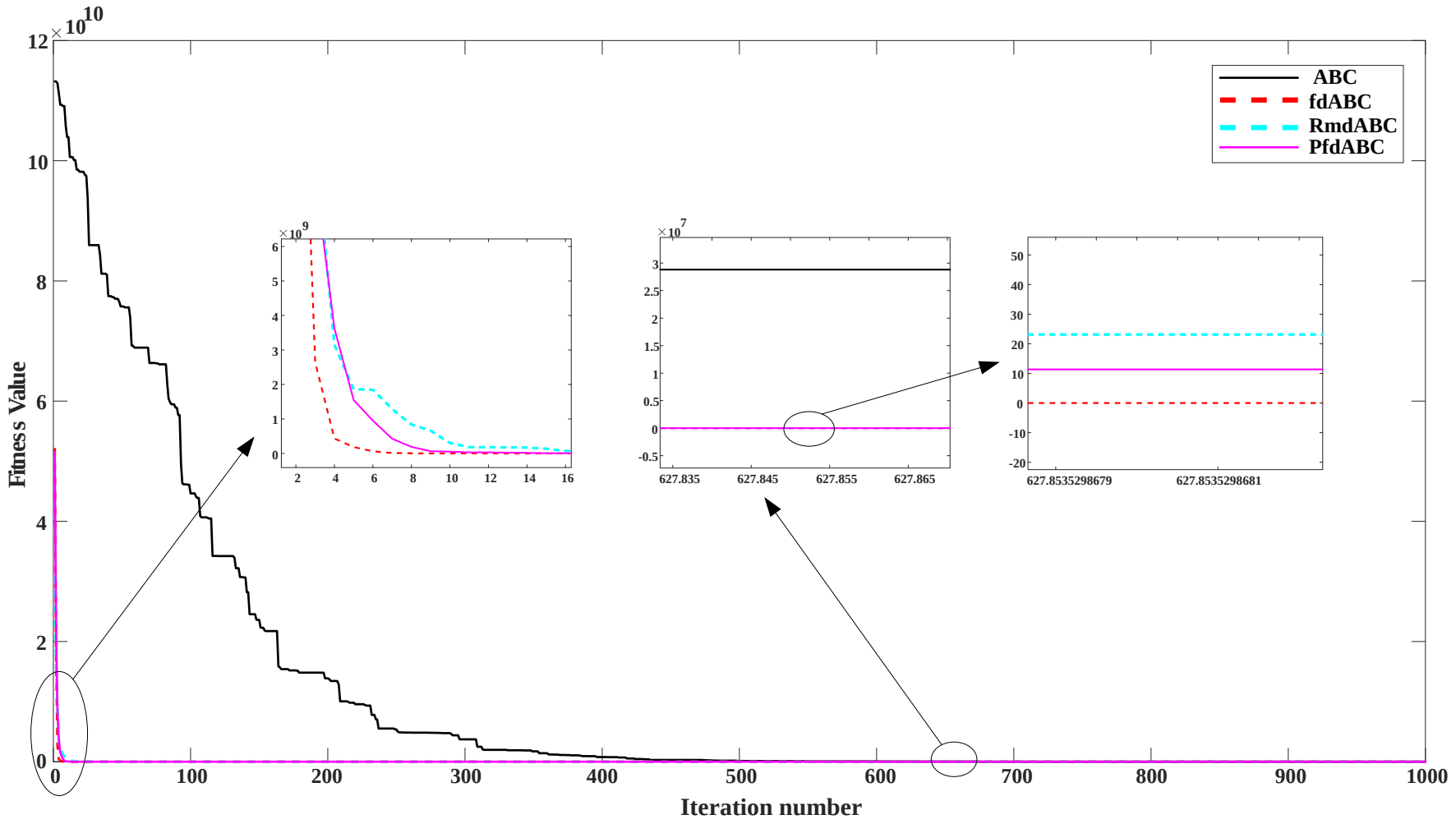

(100)

**Fig.2.** The Fitness Values with Different Strategies for Expanded Schaffer's Function.

# 4. Tasks schedule for air cargo terminal with improved ABC algorithm

The proposed ABC algorithms have proved their performance on solving complex optimization problems, in this section, they are applied to optimize the sequence of task set in air cargo terminal.

There are 16 entrances and exits totally in the container storage area, the coordinates of entrances are $R_1$(1-1-5), $R_2$(2-1-15), $R_3$(1-1-20), $R_4$(1-1-25), $R_5$(1-1-30), $R_6$(1-1-35), $R_7$(1-1-40), $R_8$(1-1-50), $R_9$(1-1-60) and the coordinates of exits are $C_1$(1-1-8), $C_2$(1-1-18), $C_3$(1-1-28), $C_4$(1-1-38), $C_5$(1-1-48), $C_6$(2-1-53), $C_7$(1-1-58). Sixty tasks are waiting for being scheduled where the first thirty tasks are inbound and the later thirty ones are outbound, and the corresponding positions to be assigned are defined as Table 7, where the first value in the bracket represents the row number of shelf; the second value indicates the number of layers and the third value indicates the number of columns.

**Table 7.** Task sets to be scheduled.

| Input tasks | | Input tasks | | Output tasks | | Output tasks | |
|---|---|---|---|---|---|---|---|
| 1 | I(1-5-34) | 16 | I(1-8-44) | 31 | O(1-3-10) | 46 | O(2-8-10) |
| 2 | I(2-3-14) | 17 | I(2-8-32) | 32 | O(1-5-55) | 47 | O(1-3-32) |
| 3 | I(1-3-58) | 18 | I(2-3-54) | 33 | O(1-5-25) | 48 | O(1-4-50) |
| 4 | I(1-5-26) | 19 | I(1-3-40) | 34 | O(2-4-8) | 49 | O(2-3-38) |
| 5 | I(1-5-30) | 20 | I(1-4-60) | 35 | O(2-2-18) | 50 | O(2-1-58) |
| 6 | I(1-2-55) | 21 | I(1-3-20) | 36 | O(2-1-16) | 51 | O(1-5-24) |
| 7 | I(1-5-24) | 22 | I(2-2-43) | 37 | O(2-3-51) | 52 | O(1-4-30) |
| 8 | I(1-4-40) | 23 | I(2-4-50) | 38 | O(1-5-6) | 53 | O(2-6-40) |
| 9 | I(1-5-40) | 24 | I(1-6-10) | 39 | O(2-5-3) | 54 | O(2-4-35) |
| 10 | I(1-5-35) | 25 | I(2-7-20) | 40 | O(1-6-12) | 55 | O(2-8-51) |
| 11 | I(2-5-23) | 26 | I(1-6-15) | 41 | O(2-6-13) | 56 | O(2-2-30) |
| 12 | I(1-7-43) | 27 | I(2-8-30) | 42 | O(2-7-49) | 57 | O(1-2-60) |
| 13 | I(1-3-48) | 28 | I(2-2-45) | 43 | O(1-7-57) | 58 | O(1-3-26) |
| 14 | I(1-8-50) | 29 | I(1-7-58) | 44 | O(1-5-25) | 59 | O(1-6-35) |

| 15 | I(1-6-21) | 30 | I(1-4-9) | 45 | O(2-6-18) | 60 | O(2-8-10) |
|---|---|---|---|---|---|---|---|

Before scheduling, the mapping relationship between the solution of the scheduling problem and the food source of ABC algorithm needs to be established, and a sort mapping coding (SMC) strategy is introduced. It assigns a random number to each dimension of the solution, then sorts the numbers in ascending order and assigns index values based on their sequence. The resulted sequence represented by the index value is the scheduling scheme. The solution generated by SMC encoding is shown in Table 8. The number of the second row is randomly generated in each dimension of a solution, and the number of the third row which is the index of random number represents the scheduling order of cargoes.

The solution represented by SMC can be updated in each iteration of algorithm. During the phases of employed bee and the onlooker bees, the value of certain dimension can be updated based on Eq.7, and the number of dimensions to be updated is determined with different ABC algorithms: ABC algorithm only selects one dimension, all dimensions can be updated in fdABC and PfdABC, RmdABC randomly selects certain dimensions.

**Table 8.** Generated scheduling scheme with SMC.

| 1 | 2 | 3 | 4 | 5 | 6 | 7 | 8 | 9 | 10 | 11 | 12 | 13 | 14 | 15 |
|---|---|---|---|---|---|---|---|---|---|---|---|---|---|---|
| -2.47 | -8.65 | 9.0 | -7.9 | -2.25 | 0.14 | 7.88 | 5.22 | -2.98 | 7.36 | 2.18 | -8.60 | -1.80 | 7.08 | 1.62 |
| 26 | 8 | 58 | 12 | 29 | 36 | 56 | 50 | 24 | 54 | 43 | 9 | 31 | 53 | 40 |
| 16 | 17 | 18 | 19 | 20 | 21 | 22 | 23 | 24 | 25 | 26 | 27 | 28 | 29 | 30 |
| -6.30 | -9.15 | -6.0 | -5.37 | -2.84 | 5.3 | 1.81 | -1.36 | 4.78 | -4.44 | -9.27 | 1.38 | 5.46 | -3.79 | 2.31 |
| 14 | 5 | 15 | 17 | 25 | 51 | 42 | 33 | 46 | 18 | 4 | 39 | 52 | 20 | 44 |
| 31 | 32 | 33 | 34 | 35 | 36 | 37 | 38 | 39 | 40 | 41 | 42 | 43 | 44 | 45 |
| -0.91 | 4.03 | -3.34 | -8.10 | -2.15 | -3.64 | 1.78 | -6.82 | 5.17 | 7.58 | 4.82 | -2.33 | 0.72 | -9.63 | -0.85 |
| 34 | 45 | 23 | 11 | 30 | 22 | 41 | 13 | 49 | 55 | 47 | 27 | 37 | 2 | 35 |
| 46 | 47 | 48 | 49 | 50 | 51 | 52 | 53 | 54 | 55 | 56 | 57 | 58 | 59 | 60 |
| 1.36 | -8.30 | -1.69 | 9.52 | 9.80 | -4.42 | -8.73 | -9.47 | -9.12 | -2.28 | -3.69 | 8.14 | 5.04 | -5.81 | -9.82 |
| 38 | 10 | 32 | 59 | 60 | 19 | 7 | 3 | 6 | 28 | 21 | 57 | 48 | 16 | 1 |

Comparative studies among four algorithms, including ABC, fdABC, PfdABC and RmdABC, are executed for the above scheduling problem. They are performed on the computer with Inter(R) Core (TM) i7-8750h CPU @2.20Ghz, 16 GB of memory. The swarm size of corresponding algorithms is set to be 200, the maximum number of local searches is 100, and the maximum number of iterations is equal to 1500.

The scheduling program ran 20 times and Tables 9-12 show the optimization results with different algorithms respectively, the second and third rows of tables are the runtime of algorithm and task execution time corresponding to the optimal solution.

Table 13 presents the optimization results with different algorithms of the first ten trials, and important indices including the average time required to execute the sequence corresponding to optimal solution in different independent trials (Avg), the shortest inbound and outbound time (Min), the longest inbound and outbound time(Max), the corresponding computational time(CPU time). The two tables show that the four algorithms all can handle the complicated cargo task set scheduling problem, and compared with ABC, the minimum/maximum/average inbound and outbound time corresponding to PfdABC and RmdABC decreased by 4% at most, there are relatively small gaps among the average inbound and outbound time of the whole task set obtained with PfdABC or RmdABC and fdABC. On the other hand, the runtime of the optimization problem with PfdABC and RmdABC reduces by 60% compared with fdABC.

**Table 9.** The scheduling results with ABC.

| 1 | 2 | 3 | 4 | 5 | 6 | 7 | 8 | 9 | 10 |
|---|---|---|---|---|---|---|---|---|---|
| 58.83691 | 58.55166 | 59.21077 | 59.00102 | 59.48015 | 60.08423 | 59.64511 | 59.93071 | 60.0780 | 59.45079 |
| 6852.899 | 6772.394 | 6806.412 | 6697.144 | 6892.792 | 6905.639 | 6748.778 | 6865.903 | 6787.75 | 6803.033 |

| 11 | 12 | 13 | 14 | 15 | 16 | 17 | 18 | 19 | 20 |
|---|---|---|---|---|---|---|---|---|---|
| 59.2933 | 59.7117 | 60.8208 | 59.75595 | 59.4666 | 60.08446 | 60.265863 | 57.03016 | 59.3823 | 59.5545 |
| 6839.283 | 6824.894 | 6877.635 | 6796.278 | 6787.26 | 6864.1491 | 6808.0102 | 6871.288 | 6839.524 | 6817.274 |

**Table 10.** The scheduling results with fdABC.

| 1 | 2 | 3 | 4 | 5 | 6 | 7 | 8 | 9 | 10 |
|---|---|---|---|---|---|---|---|---|---|
| 3633.212 | 3400.792 | 3406.109 | 3431.726 | 3415.695 | 3369.328 | 3373.345 | 3391.797 | 3567.340 | 3492.05 |
| 6611.737 | 6613.612 | 6615.487 | 6617.482 | 6613.612 | 6619.357 | 6609.982 | 6619.116 | 6613.612 | 6615.48 |

| 11 | 12 | 13 | 14 | 15 | 16 | 17 | 18 | 19 | 20 |
|---|---|---|---|---|---|---|---|---|---|
| 3462.246 | 3458.2017 | 3472.602 | 3512.458 | 3432.370 | 3490.9916 | 3562.685 | 3492.075 | 3414.187 | 3515.0619 |
| 6617.482 | 6613.6120 | 6617.241 | 6609.982 | 6619.116 | 6613.6120 | 6615.487 | 6617.482 | 6613.612 | 6617.2417 |

**Table 2.** The scheduling results with PfdABC.

| 1 | 2 | 3 | 4 | 5 | 6 | 7 | 8 | 9 | 10 |
|---|---|---|---|---|---|---|---|---|---|
| 1063.719 | 993.6694 | 1057.934 | 987.2086 | 978.4468 | 1106.609 | 1066.315 | 1053.588 | 1019.402 | 1025.188 |
| 6621.112 | 6623.107 | 6613.612 | 6609.982 | 6617.362 | 6617.241 | 6617.362 | 6611.737 | 6617.121 | 6615.487 |

| 11 | 12 | 13 | 14 | 15 | 16 | 17 | 18 | 19 | 20 |
|---|---|---|---|---|---|---|---|---|---|
| 1056.066 | 1075.639 | 1054.186 | 1030.231 | 1107.457 | 1048.584 | 1030.156 | 1057.354 | 1005.5180 | 1017.930 |
| 6606.352 | 6617.482 | 6630.487 | 6615.366 | 6611.857 | 6617.362 | 6617.482 | 6619.237 | 6615.3667 | 6611.737 |

**Table 3.** The scheduling results with RmdABC.

| 1 | 2 | 3 | 4 | 5 | 6 | 7 | 8 | 9 | 10 |
|---|---|---|---|---|---|---|---|---|---|
| 1786.48 | 1729.01 | 1729.94 | 1679.671466 | 1663.27 | 1658.65 | 1658.49 | 1681.33146 | 1666.64 | 1658.39 |
| 6617.12 | 6617.12 | 6611.85 | 6617.121340 | 6611.85 | 6611.85 | 6615.60 | 6617.12134 | 6617.24 | 6617.36 |

| 11 | 12 | 13 | 14 | 15 | 16 | 17 | 18 | 19 | 20 |
|---|---|---|---|---|---|---|---|---|---|
| 1660.59 | 1664.006313 | 1658.16 | 1685.00 | 1691.47 | 1662.51 | 1662.37832 | 1687.67 | 1774.40 | 1782.15 |
| 6621.11 | 6611.857451 | 6617.12 | 6611.85 | 6615.60 | 6615.60 | 6615.60745 | 6611.85 | 6617.12 | 6611.85 |

**Table 4.** The scheduling results.

| | Min(s) | Max(s) | Avg(s) | CPU(s) |
|---|---|---|---|---|
| ABC | 6697.144489 | 6905.639859 | 6822.917793 | 59.4817733 |
| fdABC | 6609.982 | 6619.357 | 6615.218 | 3464.714 |
| PfdABC | 6606.352822 | 6630.487081 | 6616.342868 | 1041.760423 |
| RmdABC | 6611.857452 | 6621.112081 | 6615.193794 | 1692.015441 |

Tables 14-17 show the resulted sequence and entrances and exits allocation corresponding to the optimal solutions with different algorithms.

A.Optimal result with PdmABC: 6606.352822s

Runtime of Algorithm: 1056.066203s

Optimal scheduling route: 15 16 46 31 24 41 55 42 54 56 29 25 45 27 44 14 35 58 9 47 43 49 19 23 32 13 28 22 37 18 6 3 57 50 20 48 8 59 17 38 33 12 52 7 11 60 51 34 39 5 1 36 10 53 21 26 4 30 40 2

**Table 5.** Entrances and exits allocation scheme with PfdABC.

| | Inbound tasks | | Outbound tasks |
|---|---|---|---|
| 1 | 15,16,24,25,17,7,11,26,4,30,2 | 1 | 46,31,41,45,44,59,38,33,60,51,34,39,40 |
| 2 | 14,12,5,1,10,21 | 2 | 55,42,35,58,52,36,53 |
| 3 | 9, | 3 | 54,56,47,43,49, |
| 4 | 29,27, | 4 | 32,48, |
| 5 | 19, | 5 | 37, |
| 6 | 23, | 6 | |
| 7 | 13,28, | 7 | 57,50, |
| 8 | 18,6,3,20,8,3 | | |

B.Optimal results with fdABC: 6609.982452s

Runtime of Algorithm: 3512.458081s

Optimal scheduling route: 16 24 25 45 39 26 38 15 31 14 1 36 35 29 8 56 49 19 32 13 37 20 3 57 50 6 18 28 22 23 48 42 52 4 46 51 59 27 2 17 44 5 55 58 7 40 33 41 60 30 11 34 12 21 10 53 9 54 47 43

**Table 6.** Entrances and exits allocation scheme with fdABC.

| | Inbound tasks | | Outbound tasks |
|---|---|---|---|
| 1 | 16,24,25,26,15,4,27,2,17,7,30,11 | 1 | 45,39,38,31,46,51,59,44,40,33,41,60,34, |
| 2 | 14,1,5,12,21,10 | 2 | 36,35,42,52,55,58,53 |
| 3 | 9, | 3 | 56,49,54,47,43 |
| 4 | 29,8 | 4 | 32,48, |
| 5 | 19, | 5 | 37, |
| 6 | 23, | 6 | |
| 7 | 13,28,22 | 7 | 57,50, |
| 8 | 20,3,6,18 | | |

C.Optimal results with RmdABC: 6611.857452s

Runtime of Algorithm: 1685.007386s

Optimal scheduling route: 15 39 59 21 52 14 11 30 27 34 38 41 40 7 26 31 46 44 12 42 9 47 19 29 8 58 1 36 54 49 48 13 22 37 20 3 57 50 6 18 28 32 23 10 24 25 17 60 2 33 45 5 53 16 4 51 55 35 43 56

**Table 7.** Entrances and exits allocation scheme with RmdABC.

| | Inbound tasks | | Outbound tasks |
|---|---|---|---|
| 1 | 15,11,30,27,7,26, 24,25,17,2,16,4 | 1 | 39,59,34,38,41,40,31,60,33,45,51, 46, 44, |
| 2 | 21,14,12,1,10,5 | 2 | 52,42,58,36,53,55,35 |
| 3 | 9, | 3 | 47,54,49,43,56 |
| 4 | 29,8, | 4 | 48,32 |
| 5 | 19, | 5 | 37, |
| 6 | 23, | 6 | |
| 7 | 13,22,28 | 7 | 57,50, |
| 8 | 20,3,6,18 | | |

D. Optimal results with ABC: 6697.144489s

Runtime of Algorithm: 59.001024s

Optimal scheduling route: 17 34 25 14 33 51 10 9 22 60 16 52 36 11 15 27 40 19 55 56 45 26 23 44 7 39 21 38 24 35 41 4 30 12 20 18 3 1 54 28 49 43 53 37 8 31 42 47 32 2 6 46 48 13 57 5 59 50 58 29

**Table 8.** Entrances and exits allocation scheme with ABC.

| | Inbound tasks | | Outbound tasks |
|---|---|---|---|
| 1 | 17,25,10,16,11,15,27,26,7,24,4,30,2, | 1 | 34,33,51,60,40,45,44,39,38,41,31,46,59 |
| 2 | 14,21,12,1,5 | 2 | 52,36,55,35,54,53,42,58, |
| 3 | 9, | 3 | 56,49,43,47, |
| 4 | 8,29 | 4 | 32,48, |
| 5 | 19, | 5 | 37, |
| 6 | 23, | 6 | |
| 7 | 22,28,13, | 7 | 57,50, |
| 8 | 20,18,3,6, | | |

The results show that the four algorithms can handle the complex scheduling problem, and compared with ABC algorithm, fdABC, RmdABC and PfdABC algorithms can greatly improve the performance of solutions, the minimum/maximum/average scheduling efficiency in twenty times increased by 4% at most. Meanwhile, the proposed RmdABC and PfdABC algorithms can improve the operational program efficiency compared with fdABC, the time drops by 69.9% and the performance of solutions are almost the same. The scheduling results prove that the proposed RmdABC and PfdABC could balance the exploration and exploitation abilities efficiently.

# 5. Conclusion

To improve the efficiency of the air cargo terminal, ABC algorithm is introduced to schedule the task set in this paper. Moreover, for increasing the diversity of the population meanwhile accelerating the convergence, improved algorithms including PfdABC, RmdABC are proposed to enhance the optimization performance. The experimental results show that the ABC, fdABC, RmdABC and PfdABC can solve the air cargo terminal task sets scheduling problem efficiently, and the proposed RmdABC and PfdABC could balance exploration and exploitation performance efficiently.

Although ABC and corresponding improved algorithms have been proved to be effective in solving the complex scheduling problem, how to appropriately choose the control parameters for different optimization issues is still a problem in ABC as with other metaheuristic algorithms. This problem has not been investigated sufficiently in the literature. Therefore, designing a general principle for tuning the control parameters of ABC can be addressed as a searching subject in future studies.

## Author's Contributions

Conceptualization, Haiquan Wang, Menghao Su. Funding acquisition, Shengjun Wen, Haiquan Wang. Investigation and methodology, Xiaobin Xu, Jianhua Wei, Haiquan Wang. Supervision, Ran Zhao. Writing of the original draft, Haiquan Wang, Yan Wang, Hongjun Li. Writing of the review and editing, Haiquan Wang, Ping Liu, Hans-Dietrich Haasis.

## Funding

This work was supported by High-end foreign expert program of Ministry of Science and Technology (No. G2021026006L); the Training Program for Young Teachers in Universities of Henan Province (No. 2020GGJS137); the Key Scientific Research Projects of Henan Province (No.22A413011); Henan Province Science and Technology R&D projects (No. 202102210135, 212102310547, 212102210080 and 222102210019); National Nature Science Foundation of China (No. U1813201); Key Program of University of Henan Province (No. 22A413011).

## Competing Interests

The authors declare that they have no competing interests.